\documentclass[12pt]{article}
\usepackage{epsfig}
\usepackage{amsmath}

\newcommand{\nit}{\noindent}
\newcommand{\no}{\nonumber}
\newcommand{\be}{\begin{equation}}
\newcommand{\ee}{\end{equation}}
\newcommand{\br}{\begin{eqnarray}}
\newcommand{\er}{\end{eqnarray}}

\newcommand{\abs}[1]{\lvert #1 \rvert}
\newcommand{\norm}[1]{\lVert #1 \rVert}

\newcommand{\normLinf}[1]{{\Vert #1 \rVert}_{\infty}}

\newcommand{\avg}[1]{\langle #1 \rangle}
\newtheorem{theo}{Theorem}[section]

\newtheorem{prop}{Proposition}[section]
\newtheorem{lem}{Lemma}[section]
\newtheorem{cor}{Corollary}[section]

\begin{document}

\title{A Variational Principle Based Study of KPP Minimal Front Speeds in Random Shears}

\author{James Nolen\thanks{Department of Mathematics, University of Texas at Austin,
Austin, TX 78712 (jnolen@math.utexas.edu).}
\and Jack Xin\thanks{Department of Mathematics and ICES (Institute of Computational
Engineering and Sciences), University of Texas at Austin, Austin, TX 78712
(jxin@math.utexas.edu).}}

\setcounter{section}{0}
\date{}
\maketitle
\thispagestyle{empty}
\begin{abstract}
Variational principle for Kolmogorov-Petrovsky-Piskunov (KPP) minimal front speeds provides
an efficient tool for statistical speed analysis, as well as a fast and accurate method for speed computation.
A variational principle based analysis is carried out
on the ensemble of KPP speeds through spatially stationary
random shear flows
inside infinite channel domains. In the regime of small
root mean square (rms) shear amplitude, the enhancement of the ensemble averaged KPP front speeds is proved to
obey the quadratic law under certain shear moment conditions. Similarly, in the large rms amplitude regime,
the enhancement follows the linear law. In particular, both laws hold for the Ornstein-Uhlenbeck process
in case of two dimensional channels. An asymptotic ensemble averaged speed formula is derived in the small rms
regime and is explicit in case of the Ornstein-Uhlenbeck process of the shear. Variational
principle based computation agrees with these analytical findings, and
allows further study on the speed enhancement distributions as well as the dependence of enhancement on
the shear covariance. Direct simulations in the small rms regime
suggest quadratic speed enhancement law for non-KPP nonlinearities.

\end{abstract}

\newpage
\setcounter{page}{1}
\section{Introduction}
\setcounter{equation}{0}
Front propagation in heterogeneous fluid flows
has been an active research topic for decades (see \cite{CW}, \cite{Kh}, \cite{KA}, \cite{Ro},
\cite{Vlad}, \cite{Xin1}, \cite{Yak} and references therein).
A fascinating phenomenon is that the large time (large scale) front speed
can be enhanced due to the presence of multiple scales in fluid flows.
Speed characterizations and enhancement laws have been studied mathematically
for various flow patterns by analysis of the proto-type models, e.g.
the reaction-diffusion-advection equations (see
\cite{B2, BH1,Const1,HPS,KR,Kh,MS1,MS2,NX1,PX,Vlad,Xin1,Xin2} and references
therein). The enhancement obeys a quadratic law in the small amplitude flow regime,
known as the Clavin-Williams relation \cite{CW}, which is proved to be true for
deterministic shear flows \cite{PX}, \cite{HPS}, \cite{NX2}, \cite{NX1} etc.
The enhancement is proved to grow linearly in large amplitudes of deterministic flows (shear and
percolating flows), see \cite{B1,Const1,HZ1,KR,NX2,NX1} etc.
However, enhancement exponent 4/3 was reported based on numerical simulation of
random Hamilton-Jacobi models (so called G-equation or KPZ model) on fronts in
weak randomly stirred array of vortices \cite{KR}. Likewise, formal renormalization group method
\cite{Yak} suggested that front speeds may grow sublinearly (slower than linear by a logarithmic
factor) in the strong random flow regime.
These two findings raised the issue as to what extent
the speed enhancement laws in the deterministic flows
are valid for random flows. Recently, one of the present authors \cite{Xin2} showed that
both the quadratic and linear laws hold almost surely
for Kolmogorov-Petrovsky-Piskunov (KPP) minimal front speeds through white in time
spatially Gaussian random shears on the plane. Yet the KPP front speeds diverge almost surely in
spatially Gaussian random shears. A powerful tool is the variational principle of KPP minimal front
speeds.

In this paper, we consider KPP front speeds
through random shear flows in channel domains $D\equiv R \times \Omega$, where $\Omega \subset R^{n-1}$,
$n \geq 2$, a bounded simply connected domain with smooth boundary.
We shall address the enhancement laws of the ensemble averaged front speeds.
The model equation is:
\be
u_t  =  \Delta_{x,y} u + B \cdot \nabla_{x,y} u + f(u), \label{e1}
\ee
where $t \in R^{+}$,
$\Delta_{x,y}$ the $n$-dimensional Laplacian, $(x,y) \in D$.
The nonlinearity $f=u(1-u)$, the so called KPP reaction.
Other nonlinearities \cite{Xin1} will be discussed later.
The vector field $B=(b(y,\omega),\bf{0})$ where $b(y,\omega)$ is a stationary continuous scalar random process in $y$,
its ensemble mean equal to zero. The zero Neumann boundary condition is imposed at
$\partial \, \Omega $: $\frac{\partial u}{\partial \nu} = 0$,
$\nu$ the unit outward normal.

For nonnegative initial data approaching zero and one at $x$ infinities rapidly enough,
the KPP solutions propagate as fronts with speed $c^{*}$ given by the variational principle \cite{BN1,Xin1,BH1}:
\be
c^* = c^*(\omega) =  \inf_{\lambda>0} \frac{\mu(\lambda, \omega)}{\lambda},  \label{e8}
\ee
where $\mu(\lambda, \omega)$ is the principal eigenvalue with corresponding eigenfunction $\phi > 0$ of the
problem:
\br
\bar L_\lambda \, \phi = \Delta_y \; \phi + [\lambda^2 + \lambda \, b(y,\omega) + f'(0)]\, \phi & = &
\mu(\lambda, \omega)\, \phi , \quad y \in \Omega, \label{e7} \\
\frac{\partial \phi}{\partial \nu} & = & 0, \quad y \in \partial
\Omega. \label{e23}
\er

The variational speed formula (\ref{e8}) makes possible an analysis of
ensemble averaged random front speeds. Using variational formulas on the principal eigenvalue $\mu$, we are able to
obtain tight upper and lower bounds on $c^{*}$ and estimate $E[c^{*}]$ in terms of
moments of suitable norms of the shear over $\Omega$. For small root mean square (rms) shear, the quadratic enhancement
law is proved and is explicit in case of the Ornstein-Uhlenbeck (O-U) process when $n=2$. The linear growth law holds
in the large rms regime under weaker moment conditions on the shear. In both regimes, the moment conditions are satisfied by
the O-U process when $n=2$.

The variational formula (\ref{e8}) also offers an efficient and accurate way to
compute a large ensemble of random front speeds.
Directly solving the original time dependent equation (\ref{e1}) to reach
steady propagating states can be both slow and less accurate. Numerical difficulties abound for direct simulations
due to three large parameters. A large ensemble of random fronts requires a large enough
truncated channel domain to contain the front over large times. Due to
occasional random excursions in $b$, the domain size in $x$ has also to be made adaptively large.
This can be prohibitively expensive in the regime of large rms shears.

The variational formula (\ref{e8}) allows us to compute quickly and accurately the ensemble averaged speeds in both
small and large shear rms regimes when $n=2$. An interesting difference from the deterministic case is that
the integral average of $b(y,\omega)$ in $y \in \Omega$, i.e. $\bar{b}=\bar{b}(\omega)=
|\Omega|^{-1}\int_{\Omega}\, b(y,\omega)\, dy$, is a random constant not equal to zero.
This quantity can be of either sign, and influence greatly the numerical approximation of $E[c^{*}]$
in the small rms regime, even though it does not contribute to the exact $E[c^{*}]$, since $E[\avg{b}]=0$.
To assess the speed enhancement accurately, we subtract this random constant from each $c^{*}(\omega)$
before evaluating the expectation numerically.
This way, we are able to minimize the errors in approximating $E[c^{*}]$ in a finite ensemble.
In our computation, $b$ is a discrete O-U process.

In complete agreement with analysis, we find numerically that the
ensemble averaged speeds obey a quadratic
law in the small rms regime and a linear law in the large rms regime.
Without the $\bar{b}$ subtraction technique, the computed average speed enhancement in the small rms regime can give
inaccurate scaling exponents significantly below two. The same technique and direct simulations for other nonlinearities
(combustion, bistable) suggest quadratic speed enhancement in the small rms regime.
The computed speed ensemble then permits us to study further the enhancement distribution and
its dependence on variation of shear covariance functions.

The paper is organized as follows. In section 2, we prove the enhancement laws of ensemble averaged speeds based on
variational principles, and derive a closed form speed asymptotic formula in case of O-U process.
In section 3, we describe numerical methods for computing speed ensemble with the variational
formula (\ref{e8}) and speed statistics, then show numerical results.
We also make comparisons with the prediction of the asymptotic formula and with
direct simulations. In section 4, we conclude with a remark on future works.

\section{Average Speed Asymptotics}
Consider scaling shear amplitude $b(y) \mapsto \delta b(y)$,
and denote by $c^*(\delta)$ the minimal KPP speed corresponding to shear $\delta b$.
Let $c_0 = c^*(0)= 2 \sqrt{f'(0)}$ denote the minimal speed in case of zero advection.
If the shear $b = b(y)$ has zero integral average over $\Omega$,
$
\avg{b} = \frac{1}{|\Omega |}\; \int_{\Omega}\, b(y)\;dy = 0,$
the corresponding minimal speed $c^{*}(\omega)$ is always enhanced by the shear. This is
true also for time dependent shears, see \cite{NX1,NX2} and references therein.
For each realization, $c^*(\delta, \omega) = c^*_0 + O(\delta^2)$ as $\delta \ll 1$;
$c^*(\delta, \omega) = c^*_0 + O(\delta)$ as $\delta \gg 1$.
For each realization, as $\delta \ll 1$, we have:

\begin{prop}
Let $\chi =\chi(y)$ solve $\Delta_y \chi =  -b$, $y \in \Omega$, with zero Neumann boundary condition, where
$b \in C(\overline{\Omega})$, has zero mean over $\Omega$.
Then for $\delta$ sufficiently small, the minimal speed has the expansion
\be
c^*(\delta) = c_0 + \frac{c_0 \delta^2}{2 \abs{\Omega}} \int_\Omega \abs{\nabla \chi}^2\;dy + O(\delta^3). \label{e16}
\ee
\end{prop}

Up to $O(\delta^2)$, the formula (\ref{e16}) is independent of the
nonlinearity, see \cite{PX,HPS} for $f$ being
bistable or combustion nonlinearity,
also \cite{NX1} for more general nonlinearities and time
periodic shears.  We shall give a proof of Proposition 2.1
using variational formulas, and later generalize it to the random case.
A helpful fact is that
the infimum in (\ref{e8}) can be restricted to a bounded set independent
of $b$ and $\delta$, as stated in the following lemma:
\begin{lem}
Let $b \in C(\overline{\Omega})$ have zero mean over $\Omega$, and let $\lambda_0 = \sqrt{f'(0)}$. Then
\be
\inf_{\lambda>0} \frac{\mu(\lambda)}{\lambda} =  \inf_{0 < \lambda \leq \lambda_0} \frac{\mu(\lambda)}{\lambda}.
\ee
\end{lem}

\begin{figure}[tb]
\centerline{\epsfig{file=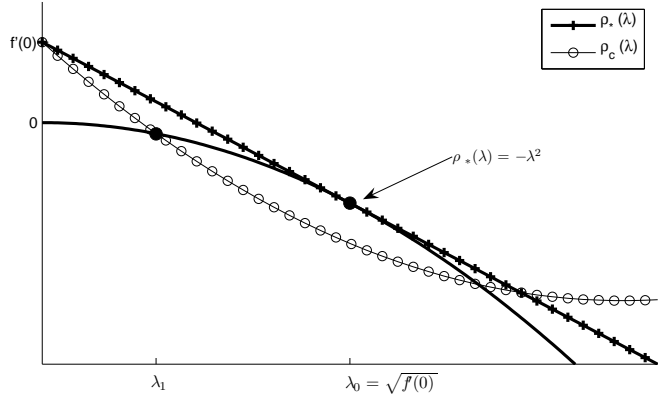,width=300pt}}
\caption{Intersecting curves $\rho_c$ and $-\lambda^2$.}
\label{fig14}
\end{figure}

\nit {\it Proof of lemma:}
For each $c>0$, we let $\rho_c(\lambda) = \mu(\lambda) - \lambda c - \lambda^2$. So, if $\phi>0$ is the eigenfunction defined by (\ref{e7}), then $\rho_c(\lambda)$ is
the principal eigenvalue defined by the equation
\br
\Delta_y \phi + [ \lambda b(y)- \lambda c + f'(0)]\phi = \rho_c(\lambda) \phi , \quad y \in \Omega .
\er
One can readily verify that $\partial_{\lambda}\;
\rho_{c}(\lambda)|_{\lambda =0} = -c < 0$.
The variational formula (\ref{e8}) can be expressed as
\br
c^* &=& \inf \left\{ c \; \lvert \; \exists \lambda > 0 ,  \lambda c = \mu(\lambda) \right\} \no \\
& = &  \inf \left\{ c \; \lvert \; \exists \lambda > 0 , \rho_c(\lambda) = -\lambda^2  \right\} \label{e32}.
\er
Consider the points where $\rho_c(\lambda) = -\lambda^2$.
By Proposition 2.1 of \cite{BN1},
the continuous curve $\lambda \mapsto \rho_c(\lambda)$ is
convex in $\lambda$, for each $c>0$.
Also, $\rho_c(0) = f'(0) > 0$. Therefore,
for a given $c>0$, there can be at most two values of $\lambda > 0$ such
that $\rho_c(\lambda) = -\lambda^2$. The line $\rho_*(\lambda) = - 2 \sqrt{f'(0)} \lambda + f'(0)$ satisfies $\rho_*(\lambda) \geq -\lambda^2$, with equality
holding only at one point: $\lambda_0 = \sqrt{f'(0)}$. Since $\rho_*(0) = \rho_c(0)$ and $\rho_c(\lambda)$ is convex and $\rho_*$ is a line,
$\rho_c(\lambda) = -\lambda^2$ for some $\lambda > 0$ only if $\rho_c(\lambda_1) = -\lambda_1^2$ for some $\lambda_1 \in (0, \lambda_0]$.  This point is illustrated in Figure \ref{fig14}.
The solid curve represents the parabola $-\lambda^2$.  If $\rho_c(\lambda)$ intersects $-\lambda^2$, then one of the intersection points must be to the left of $\lambda_0 = \sqrt{f'(0)}$.
Therefore, from (\ref{e32}),
\br
\left\{ c \; \lvert \; \exists \lambda > 0 , \rho_c(\lambda) = -\lambda^2  \right\} = \left\{ c \; \lvert \; \exists \lambda \in (0, \lambda_0] , \rho_c(\lambda) = -\lambda^2  \right\} .
\er
So, we conclude that
\be
c^*(\delta) =  \inf_{0 < \lambda} \frac{\mu(\lambda)}{\lambda} =  \inf_{0 < \lambda \leq \lambda_0} \frac{\mu(\lambda)}{\lambda}. \label{e18}
\ee

\nit {\it Proof of Proposition 2.1:}

To estimate $c^*(\delta)$, we bound the principal eigenvalue $\mu(\lambda)$
using two different representations of $\mu$. First, since $\bar L$ is a self-adjoint operator, we have:
\be
\mu = \sup \frac{(\bar L_\lambda \psi,\psi)}{\norm{\psi}^2_2}, \label{e19}
\ee
where the supremum is taken over all $\psi \in H^2(\Omega)$ such that $\frac{\partial \psi}{\partial \nu} = 0$ on $\partial \Omega$.
The other representation is:
\be
\mu = \inf_{\psi} \sup_{y\in \Omega} \frac{\bar L_\lambda \psi}{\psi} = \inf_{\psi} \sup_{y\in \Omega} \left( \frac{\Delta \psi}{\psi} + \lambda \delta b  + \lambda^2 + f'(0)\right), \label{e20}
\ee
where the infimum can be taken over all $\psi \in C^1(\Omega)$ such that $\Delta \psi \in C(\Omega)$, $\psi > 0$, and
$\frac{\partial \psi}{\partial \nu} = 0$ on $\partial \Omega$.
This representation follows from the fact that the eigenfunction $\phi>0$ lies in the kernel of the self-adjoint operator $(\bar L_\lambda - \mu(\lambda) I) = (\bar L_\lambda - \mu(\lambda) I)^*$. So, if we have the strict inequality
\be
\bar L_\lambda \psi - \mu(\lambda) \psi = m < 0
\ee
then the Fredholm alternative implies that $(\phi, m)_{L^2} = 0$, a contradiction since $\phi > 0$, $m<0$. Hence,
\be
\sup_{y\in \Omega} \frac{\bar L_\lambda \psi}{\psi} \geq \mu(\lambda).
\ee
Since $\bar L_\lambda \phi = \mu(\lambda)\phi$, the formula (\ref{e20}) follows. Note that we do not require the test functions $\psi$ to be $C^2(\Omega)$, only $\Delta \psi \in C(\Omega)$. This is important since we do not want to require the shear $b(y)$ to be any more regular that $b \in C(\bar \Omega)$.

Let us derive upper and lower bounds for $\mu(\lambda)$ by choosing test functions $\psi$ as:
\be
\psi = 1 + \lambda \delta \chi + \lambda^2 \delta^2 h , \label{e21}
\ee
where $\chi = \chi(y) $ and $h = h(y)$ solve
\br
\Delta \chi &=&  -b \no \\
\Delta h &=& -b \chi + k , \label{e22}
\er
with zero Neumann boundary conditions at $\partial \, \Omega$, and $k$ a constant equal to:
\be
k = \frac{1}{\abs{\Omega}} \int_{\Omega} b\, \chi \;dy = \frac{1}{\abs{\Omega}}
\int_{\Omega} \abs{\nabla \chi}^2 \;dy . \label{e22a}
\ee
We normalize $\chi $ and $h$ so that:
\be
\inf_{x \in \Omega} \chi(x) = 0 \quad \text{and} \quad  \inf_{x \in \Omega} h(x) = 0. \label{norm1}
\ee
Then
\be
\bar L_\lambda \psi = \lambda^2 \delta^2 k + \lambda^3 \delta^3 b\,
 h + (\lambda^2 + f'(0)) \psi , \no
\ee
and
\be
\frac{(\bar L_\lambda \psi,\psi)}{\norm{\psi}^2_2} = \lambda^2 \delta^2 k \frac{\int \psi}{\int \psi^2} + \lambda^3 \delta^3 \frac{\int b \, h \psi}{\int \psi^2} + \lambda^2 + f'(0) . \label{e24}
\ee
Using the definition of $\psi$, we see that $\psi = \psi^2 - \lambda \delta \chi \psi - \lambda ^2 \delta^2 h \psi $, and
\br
\frac{\int_{\Omega} \psi \; dy}{\int_{\Omega}  \psi^2 \; dy} &=& 1 - \lambda \delta \frac{\int_{\Omega} \chi \psi \; dy}{\int_{\Omega} \psi^2 \; dy} - \lambda^2 \delta^2 \frac{\int_{\Omega} h \psi \; dy}{\int_{\Omega} \psi^2 \; dy}. \no
\er
Now from (\ref{e19}) and (\ref{e24}) we have the lower bound
\be
\mu(\lambda) \geq  \lambda^2 + f'(0)  + \lambda^2 \delta^2 k + R_1,
\ee
with
\be
R_1 = - \frac{\lambda^3 \delta^3}{\int_\Omega \psi^2} \left( k \int_{\Omega} \chi \psi +
k \lambda \delta \int_{\Omega} h \psi - \int_{\Omega} b h \psi \right) \label{e28}
\ee
By choice of $\chi \geq 0$ and $h \geq 0$, we have $\int_\Omega \psi^2 \geq \abs{\Omega}$ for all $\delta \geq 0$ and $\lambda >0$. Hence, $R_1 = O(\delta^3)$ for $\lambda$ bounded.
Returning to the variational formula (\ref{e18}), we now have a lower bound on $c^*(\delta)$:
\br
c^*(\delta) = \inf_{0 < \lambda \leq \lambda_0} \frac{\mu(\lambda)}{\lambda}
&\geq& \inf_{0 < \lambda \leq \lambda_0}
\left( \lambda + \frac{f'(0)}{\lambda}  + \lambda \delta^2 k + {R_1 \over \lambda}  \right) \no \\
&\geq & \inf_{\lambda > 0} \left( \lambda + \frac{f'(0)}{\lambda}  + \lambda \delta^2 k  \right)  + O(\delta^3)  \no \\
& = & 2 \sqrt{ f'(0) (1 + \delta^2 k)} + O(\delta^3) \no \\
&=& c_0  +  \frac{c_0 \delta^2 k}{2} + O(\delta^3). \label{e25}
\er

To obtain an upper bound on $c^*(\delta)$, we use (\ref{e20}) and calculate:
\br
\frac{\bar L_\lambda \psi}{\psi} &=&  \frac{\Delta \psi}{\psi} + \lambda \delta b  + \lambda^2 + f'(0) \no \\
& = & \frac{\lambda^2 \delta^2 k + \lambda^3 \delta^3 b\, h}{1 + \lambda \delta \chi + \lambda^2 \delta^2 h} + \lambda^2 + f'(0). \label{e29a}
\er
Since $\chi \geq 0$ and $h \geq 0$, we see from (\ref{e20}) and (\ref{e29a}) that
\be
\mu(\lambda) \leq  \sup_{y\in \Omega} \frac{\bar L_\lambda \psi}{\psi} \leq  \lambda^2 + f'(0) + \lambda^2 \delta^2 k + R_2,
\ee
with
\be
R_2 = \lambda^3 \delta^3 \norm{b h}_\infty. \label{e29}
\ee
The variational formula (\ref{e18}) implies:
\br
c^*(\delta) = \inf_{0 < \lambda \leq \lambda_0} \frac{\mu(\lambda)}{\lambda}  &\leq& \inf_{0 < \lambda \leq \lambda_0}\left( \lambda + \frac{f'(0)}{\lambda}  + \lambda \delta^2 k + R_2 \right) \no \\
& = &  c_0 +  \frac{c_0 \delta^2 k}{2} + O(\delta^3), \label{e33}
\er
completing the proof.

When the shear $b(y,\omega)$ is a random process, the corresponding minimal speed $c^*(\delta) = c^*(\delta,\omega)$ is a random variable for each $\delta$, and
we consider how the expectation $E[c^*(\delta)]$ scales with the parameter $\delta$ by
finding an exponent $p$ such that $E[c^*(\delta)] = c^*(0) + O(\delta^p)$.
Each realization of the process $b(y,\omega)$ restricted to the domain $\Omega$ does not necessarily have
zero integral over $\Omega$. Nevertheless, each realization can be written in the form
\be
b(y,\omega) = \bar b(\omega) + b_1(y,\omega), \label{e13}
\ee
where $\bar b(\omega) = \avg{b(y,\omega)}$ is the mean of $b$ over $\Omega$, and $b_1(y,\omega)$ is the variation about the mean value.
For a fixed realization, the minimal speed $c^*(\delta)$ will be affected by both the scaling of the mean $\bar b(\omega)$ and the scaling of the
varation $b_1(y,\omega)$. That is,
\br
c^*(\delta,\omega) &=& c^*_0 - \delta \bar b(\omega) + M(\delta,\omega), \label{e14}
\er
where the remainder $M(\delta,\omega)$ is the enhancement due to the variation $b_1(y,\omega)$, different for each realization. Taking the expectation
of both sides of (\ref{e14}), we have
\be
E[c^*(\delta)] = c^*_0 - \delta E[\bar b(\omega)] + E[M(\delta,\omega)] \label{e15}
\ee
For each sample, $M(\delta,\omega)$ is $O(\delta^2)$ for $\delta$ small.
Though $E[M(\delta)]$ might exhibit different scaling than quadratic, we show that the quadratic scaling law remains
for enhancement of averaged front speeds under suitable moment conditions of the shear.

\begin{theo}
Let $b(y,\omega)$ be a stationary random process in $R^{n-1}$
($n \geq 2$) so that sample paths are almost surely continuous; and that
\be
E[\norm{b}_{\infty}^{6}] < +\infty. \label{con1}
\ee
Then for $\delta$ small, the expectation $E[c^*(\delta)]$ has the expansion
\be
E[c^*(\delta)] = c_0 - \delta E[\avg{b}] + \frac{c_0 \delta^2}{2 \abs{\Omega}}
\int_\Omega E[\abs{\nabla \chi}^2]\;dy + O(\delta^3), \label{e36b}
\ee
where $b(y,\omega) = \avg{b}(\omega) + b_1(y,\omega)$;
and $\chi =\chi(y,\omega)$ solves $\Delta_y \chi =  -b_1, \; y \in \Omega$,
subject to zero Neumann boundary condition.
\end{theo}

\nit {\it Proof:}
As the contribution of $\avg{b}$ to $c^*$ is just an additive constant, it suffices to consider shear flow $b_1$
and show that it gives the averaged speed
\be
E[c^*(\delta)]  =  c_0 + \frac{c_0 \delta^2}{2 \abs{\Omega}} \int_\Omega E[\abs{\nabla \chi}^2]\;dy + O(\delta^3).
\label{e36a}
\ee
We adapt the proof of Proposition 2.1, noting however that in the stochastic case the remainders $R_1$ and
$R_2$ defined by (\ref{e28}) and (\ref{e29}) are random and not bounded
uniformly for all realizations. Instead, we will show that for $\lambda$ in a bounded interval,
\[
E[\abs{R_1}] \leq O( \delta^3) \quad \text{and} \quad E[\abs{R_2}] \leq O(\delta^3).
\]
 To this end, we estimate $\chi $ and $h$, with $C$ denoting a generic positive constant
depending only on the domain $\Omega$ and its dimension. Let $\chi$ and $h$ solve (\ref{e22}) with $\avg{\chi} = \avg{h} = 0$. Applying $W^{2,p}$ estimates (e.g. Thm. 19.1 of \cite{Fr}), we have
\be
\norm{\chi}_{W^{2,p}(\Omega)} \leq C \norm{b_1}_{L^p(\Omega)} \leq C \abs{\Omega}^{1/p} \norm{b_1}_\infty
\ee
and
\br
\norm{h}_{W^{2,p}(\Omega)} & \leq & C \norm{b_1 \chi + k}_{L^p(\Omega)} \no \\
& \leq & C \norm{b_1}_\infty \norm{ \chi }_{L^p(\Omega)} + C k \abs{\Omega}^{1/p} \no \\
& \leq & C \norm{b_1}^2_\infty
\er
since
\be
k = \avg{\abs{\nabla \chi}^2} \leq C \norm{b_1}^2_\infty
\ee
Given $\alpha \in (0,1)$, we can choose $p>1$ sufficiently large, depending on $n$, such that $W^{2,p}(\Omega)$ embeds continuously into $C^{1,\alpha}(\bar \Omega)$. It follows that there is a constant $C>0$ independent of $b$ such that
\br
\norm{\chi}_{C^1(\bar \Omega)} \leq  C \norm{b_1}_\infty \quad \text{and} \quad \norm{h}_{C^1(\bar \Omega)} \leq C \norm{b_1}^2_\infty \label{e37}
\er
If instead we normalize $\chi$ and $h$ by (\ref{norm1}), then the bounds (\ref{e37}) still hold, with different constants. Note that adding a constant to $\chi$ and $h$ does not alter the quantity $\int_\Omega \abs{\nabla \chi}^2$ that appears in the asymptotic expansion.

Now by (\ref{e37}), the integrals in $R_1$ are easily bounded as:
\br
\int_\Omega \chi \psi & = & \int_\Omega \chi + \lambda \delta \chi^2 + \lambda^2 \delta^2 \chi h \no \\
& \leq & C \left(\normLinf{b_1} +
\lambda \delta \normLinf{b_1}^2 + \lambda^2 \delta^2 \normLinf{b_1}^3
\right). \no
\er
Similarly,
\br
\int_\Omega h \psi & = & \int_\Omega h + \lambda \delta \chi h + \lambda^2 \delta^2 h^2 \no \\
& \leq & C \left(\normLinf{b_1}^2 + \lambda
\delta  \normLinf{b_1}^3
+ \lambda^2 \delta^2 \normLinf{b_1}^4 \right), \no
\er
and
\br
\left \lvert \int_\Omega b_1 h \psi \;dy \right \rvert & = & \left \lvert \int b_1 h + \lambda \delta \int b_1 h \chi + \lambda^2 \delta^2 \int b_1 h^2 \right \rvert \no \\
& \leq & C \left(\normLinf{b_1}^3 + \lambda \delta \normLinf{b_1}^4  + \lambda^2 \delta^2 \normLinf{b_1}^5 \right) . \no
\er

Since $\chi$ and $h$ are nonnegative, $\int_\Omega \; \psi^2 \; dy \geq C \int_\Omega \psi = C \abs{\Omega} > 0$, for any realization. So for $\lambda$ in a bounded interval and $\delta $ small, we bound (\ref{e28}) by:
\[ |R_1 | \leq C\, \delta^3 \, \lambda^3  \, ( 1 + \|b_1 \|_{\infty}^6 ), \]
so that
$
E[\abs{R_1}] \leq O(\delta^3).
$

To bound $E[\abs{R_2}]$, we use the normalization (\ref{norm1}) and the above estimates:
\br
\abs{R_2} & = & \lambda^3 \, \delta^3\, \norm{b_1 h}_\infty \leq C\,
\lambda^3 \, \delta^3 \, \norm{b_1}^3_\infty,
 \label{e31a}
\er
Hence $E[R_2] \leq O(\delta^3)$ for $\lambda $ in finite interval.  Now we return to (\ref{e25}) to conclude
\br
E[c^*(\delta)] &\geq& E[2\sqrt{f'(0)(1 + \delta^2 k)}] + O(\delta^3) \no \\
& = & c_0 + \frac{c_0 \delta^2 E[k]}{2} + O(\delta^3) \no
\er
since $E[k^2] \leq C E[\norm{b}^4_\infty] < \infty$,

The opposite inequality follows from (\ref{e33}) since $E[R_2] = O(\delta^3) $ for $\lambda \in (0,\lambda_0)$.
Thus formula (\ref{e36a}) holds. For general $b$, not necessarily mean zero,
\br
E[c^*(\delta)] &=& c_0 + \delta E[\avg{b}] + \frac{c_0 \delta^2}{2} E[k] + O(\delta^3) \no \\
& = & c_0 + \delta E[\avg{b}] + \frac{c_0 \delta^2}{2 \abs{\Omega}} \int_\Omega E[\abs{\nabla \chi}^2]\;dy + O(\delta^3). \no
\er
The proof is finished.

In our numerical computation on front speeds in two dimensional channels,
we use the Ornstein-Uhlenbeck (O-U) process for shear $b$, so $E[\avg{b}] = 0$.
Let us show below that the O-U process, denoted by $X(y,\omega)$,
satisfies the conditions in Theorem 2.1, and so
$E[c^*(\delta)]$ scales quadratically with $\delta$ for $\delta$ small.

\begin{cor}[Explicit Average Speed Formula]
Consider the O-U

\noindent process $b(y,\omega)$ as solution of the Ito equation:
\be
dX(y) = -a\, X(y)\,  dy + \, r \, dW(y) , \;\;\; y \in [0,L], \label{e10}
\ee
where $W(y,\omega)$ is the standard Wiener process,
$X(0,\omega) = X_0(\omega)$ is a
Gaussian random variable with mean zero, and variance $\rho = r^2/(2a)$.
Then $X(y,\omega)$ satisfies the moment conditions in Theorem 2.1.
The averaged KPP front speed in the
channel $R\times [0,L]$ is given by
\be
E[c^*(\delta)] = c_0 + {c_0 \delta^2 \over 2}{\rm enh} + O(\delta^3),\;\;\; \delta \ll 1, \label{e36}
\ee
where:
\[ {\rm enh} = \frac{r^2}{2a}\, \left(e^{-aL} \left( \frac{4}{L^2 a^4} -
\frac{1}{3a^2}  \right)  + \frac{L}{3a} - \frac{4}{L^2 a^4} - \frac{5}{3a^2} +
\frac{4}{L a^3}\right).
\]
\end{cor}

\nit {\it Proof:} The O-U process is stationary and Markov.
Its sample paths are almost surely H\"older continuous though nowhere differentible.
The process can be written as
\br
b(y,\omega) = e^{-ay} b(0,\omega) + r \int_0^y e^{-a(y-s)} dW_s(\omega) ,  \label{e11}
\er
The covariance function of this process is $\rho e^{-a|y-s|}$.
Letting $g(y,\omega)$ denote the process
\be
g(y) = e^{ay} b(y,\omega) = g(0,\omega) + r \int_0^y e^{as} dW_s(\omega),  \label{e34}
\ee
we see that $g(y,\omega)$ is a martingale \cite{KS}. By Doob's martingale moment inequality \cite{KS},
for any $p \in (1,+\infty)$,
\be
E[ \sup_{0 < y < L} \abs{g(y)}^p ] \leq \left( \frac{p}{p-1} \right)^p E[\abs{g(L)}^p].  \label{e35}
\ee
Since the process $b(y,\omega)$ is Gaussian,
(\ref{e34}) and (\ref{e35}) imply that
\be
E[\norm{b}_{\infty}^{6}] \leq C \, E[|b(L)|^6] < +\infty.
\ee
Formula (\ref{e36}) now applies to the average speed. Notice that
\be
(\chi_x(x) )^2 = \int_0^x \int_0^x b_1(s)b_1(y) \;ds \;dy , \no
\ee
and
\be
E[(\chi_x(x) )^2] = \int_0^x \int_0^x E[b_1(s)b_1(y)] \;ds \;dy \label{e27}
\ee

Let us calculate $E[b_1(s)b_1(y)]$ in terms of $E[b(s)b(y)]$.
Define:
\[
g(y) = \avg{f( \cdot , y)} \quad \text{or} \quad g(s) = \avg{f(s, \cdot)}
\]
so that
\[
E[b_1(y) b_1(s)] = E[b(s)b(y)] - E[b(s)\bar b] - E[b(y) \bar b] + E[{\bar b}^2] ,
\]
\[
E[b(s) \bar b] = \frac{1}{L} \int_0^L E[b(y) b(s)] \;dy = g(s) ,
\]
\[
E[{\bar b}^2] = \frac{1}{L^2} \int_0^L \int_0^L E[b(s)b(y)]\;dy\;ds  = \avg{g} .
\]
Thus
\[
E[b_1(y) b_1(s)] = f(s,y) +\avg{g} - g(y) - g(s) .
\]
Now, we have
\br
E[(\chi_x(x) )^2] &=& \int_0^x \int_0^x E[b_1(s)b_1(y)] \;ds \;dy \no \\
&=& \int_0^x \int_0^x f(s,y) +\avg{g} - g(y) - g(s)\;ds\;dy \no \\
&=&  x^2 \avg{g} - 2x^2 \avg{g}_x + \int_0^x \int_0^x f(s,y)\;ds\;dy, \no
\er
where $\avg{g}_x$ denotes the average of $g$ over the interval $[0,x]$, for $0 < x \leq L$. Consequently, we have
\br
E[\avg{\abs{\chi_x}^2}] &=& \frac{1}{L} \int_0^L \left( x^2 \avg{g} - 2x^2 \avg{g}_x +
\int_0^x \int_0^x f(s,y)\;ds\;dy \right) \; dx \no
\er

Using the O-U covariance function, we proceed as:
\br
g(s) &=& \frac{1}{L} \int_0^L f(s,y) \;dy  = \frac{r^2}{2a} \frac{1}{L} \int_0^L e^{-a \abs{y-s}}\;dy \no \\
& = & \frac{r^2}{2a} \left( \frac{1-e^{-as}}{La} + \frac{1 - e^{-a(L-s)}}{La}\right), \no
\er
and
\br
\avg{g}_x &=& \frac{r^2}{2a} \frac{1}{x} \int_0^x ( \frac{1-e^{-as}}{La} + \frac{1 - e^{-a(L-s)}}{La} )\;ds \no \\
&=& \frac{r^2}{2a} \left(\frac{2}{La} + \frac{1}{x L a^2}(e^{-ax} - 1) + \frac{1}{x L a^2}(e^{-aL} - e^{-a(L-x)})
\right). \no
\er
Letting $x=L$, we have
\[
\avg{g} =
\frac{r^2}{2a}\left(\frac{2}{La} + \frac{2}{L^2 a^2}(e^{-aL} - 1)\right).
\]
Similarly,
\[
\int_0^x \int_0^x f(s,y)\;ds\;dy  = \frac{r^2}{2a}\left(\frac{2x}{a} + \frac{2}{a^2}(e^{-ax} - 1)\right).
\]
Combining the above, we have
\br
E[\avg{\abs{\chi_x}^2}] &=& \frac{r^2}{2a} \left(\frac{2L}{3a} + \frac{2}{3 a^2}(e^{-aL} - 1) \right)\no \\
& & - \frac{r^2}{2a} \frac{2}{L} \int_0^L  \frac{2x^2}{La} + \frac{x}{ L a^2}(e^{-ax} - 1) +
\frac{x}{L a^2}(e^{-aL} - e^{-a(L-x)}) \;dx  \no \\
& & + \frac{r^2}{2a} \frac{1}{L} \int_0^L \frac{2x}{a} + \frac{2}{a^2}(e^{-ax} - 1)\;dx  \no \\
&=& \frac{r^2}{2a} \left(e^{-aL}
\left( \frac{4}{L^2 a^4} - \frac{1}{3a^2} \right) +
\frac{L}{3a} - \frac{4}{L^2 a^4} - \frac{5}{3a^2} + \frac{4}{L a^3}\right).
\no \\
 \label{e26}
\er
In view of (\ref{e36}), the proof is complete.

\begin{theo}[Linear Growth]
If the stationary shear process $b(y,\omega)$ has almost surely continuous sample paths and
satisfies $E[\|b\|_{\infty}] < \infty$, then
the amplified shear field $\delta \, b(y,\omega)$ generates the
average front speed:
\[ E[|c^*(\delta,\omega)|]= O(\delta), \;\;\; \delta \gg 1. \]
Moreover, $\lim_{\delta \to \infty}\, E[|c^*(\delta,\omega)|]/\delta$ exists.
\end{theo}

\nit {\it Proof:} By Theorem 5.1 of \cite{B1}, ${|c{*}(\delta,\omega)| \over \delta} \to d^{*}(\omega)$ as
$\delta \to \infty$, where $d^*$ is finite for each $\omega$. Now recall the upper bound
$|c^{*}(\delta,\omega)| \leq |c_0|+\delta \|b \|_{\infty}$. Hence for
$\delta > |c_0|$, ${|c^{*}(\delta,\omega)| \over \delta}\leq 1 + \|b \|_{\infty} \equiv Y$,
and $E(Y) < \infty$. The dominated
convergence theorem implies that:
\[ E [ {|c{*}(\delta,\omega)| \over \delta}] \to E[d^{*}(\omega)] \leq E(Y). \]
The proof is finished. Clearly, the O-U process
satisfies the required condition for linear average speed growth.

\section{Computation by Variational Principle}
\setcounter{equation}{0}
\subsection{Numerical Methods}
\begin{figure}[tb]
\centerline{\epsfig{file=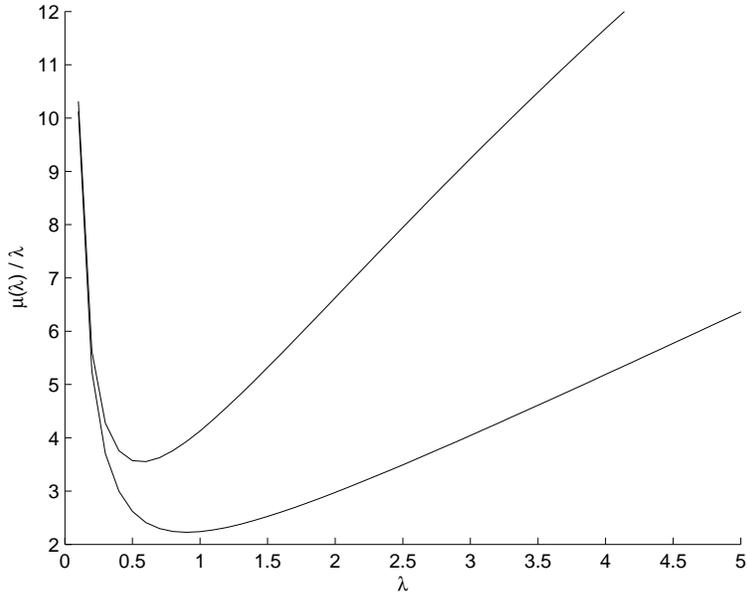,width=330pt}}
\caption{Two curves $\lambda \mapsto \frac{\mu(\lambda)}{\lambda}$}
\label{fig1}
\end{figure}

Let $n=2$. For a given $\lambda >0$, we compute the principal eigenvalue
$\mu(\lambda)$ with corresponding eigenfucntion $\phi=\phi(y) > 0$,
$y \in [0,L]$, by solving:
\br
\phi_{yy} + [\lambda^2 + \lambda \, b(y) + f'(0)]\phi = \mu(\lambda) \phi, \quad y \in (0,L),  \no \\
\frac{\partial \phi}{\partial y} = 0, \quad y=0,\, L,
\er
using a standard second order finite-difference method.
Here we suppress the random parameter $\omega$,
as computation is done realization by realization. Denote the uniform partition of the
domain by points $\left\{ y_i \right\}^m_{i=1} $, and the numerical solution by
$\bar \phi = \left\{ \bar \phi_i \right\}^m_{i=1} $, where
$h = L/(m-1)$, $y_i = (i-1)h$, and $\bar \phi_i \approx \phi(y_i)$.
The discretized system is
\be
\frac{1}{h^2} \bar \phi_{i-1} + (\lambda^2 + \lambda b_i + f'(0) -
\frac{2}{h^2})\bar \phi_i + \frac{1}{h^2}\bar \phi_{i+1}
= \mu(\lambda)\bar \phi_i \quad i = 2,\dots,m-1 , \no
\ee
with second order approximation of the Neumann boundary conditions.
This reduces to finding the principal eigenvalue of a symmetric tridiagonal matrix,
easily accomplished with double precision LAPACK routines \cite{LAP}.
Then we compute points on the
curve $H(\lambda) = \frac{\mu(\lambda)}{\lambda}$, and
minimize over $\lambda$ using a Newton's method with line search.
Our approximation decreases with each iteration and
converges quadratically in the region near the
infimum. Two illustrative curves $H(\lambda)$ are shown in Figure \ref{fig1}
for two different realizations of the shear.

\begin{figure}[tb]
\centerline{\epsfig{file=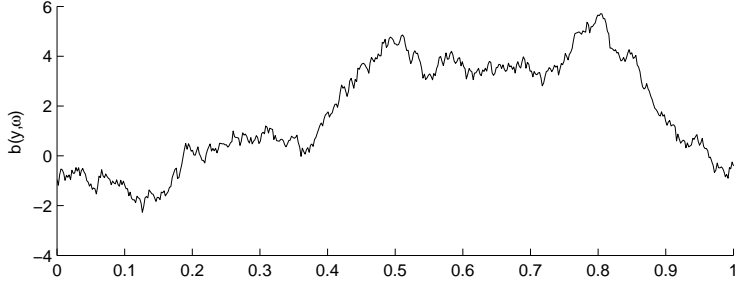,width=330pt}}
\caption{One sample path of the Ornstein-Uhlenbeck process $b(y,\omega)$.}
\label{fig7}
\end{figure}

\begin{figure}[tb]
\centerline{\epsfig{file=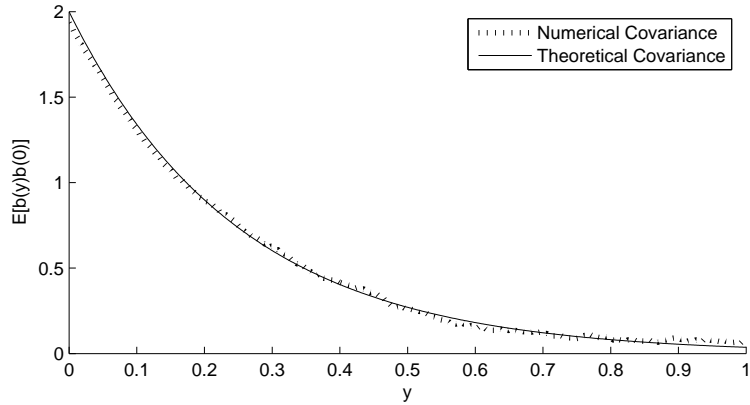,width=330pt}}
\caption{Numerical and exact covariance functions of the Ornstein-Uhlenbeck process $b(y,\omega)$.}
\label{fig6}
\end{figure}

We generate realizations of the shear process $b(y,\omega)$  by numerically evaluating the stochastic ODE (\ref{e10}) with the Milstein
scheme (see \cite{KoPl}). Although this scheme is first order, we use a discrete spacing $\bar h \leq h^2$, where $h$ is the discrete grid spacing for
the eigenvalue problem, so that the method is still
second order accurate in the parameter $h$.
Figure \ref{fig7} shows a sample path, and Figure \ref{fig6} compares exact
and numerical covariance function constructed from 5000 samples.

To approximate the expectation $E[c^*(\delta)]$ we generate $N$ independent realizations (indexed by $i=1,\dots,N$) of the shear and
compute the corresponding minimal speeds $\{ c^*_i \}$ for each $\delta$. Then we compute the average
\be
E[c^*(\delta)] \approx \bar E(\delta) =  c^*_0 + \frac{1}{N} \sum^N_{i=1} M_i(\delta),
\ee
where $M_i(\delta) = c^*_i(\delta) - c^*_0 + \delta \bar b_i$. That is, we subtract the linear part due to the mean of the shear being nonzero, as in (\ref{e14}).

Once we have the averages $\bar E(\delta)$ for each $\delta$, we compute the exponents $p$
using the least squares method to fit a line to a log-log plot
of speed versus amplitude.  That is, the exponent $p$ is the slope of the best-fit line
through the data points $(\log(\delta), \log(\bar E[c^*(\delta)] - c^*(0)))$ for
each shear amplidute $\delta$.

\subsection{Numerical Results}
\setcounter{equation}{0}

\subsubsection{Scaling with shear amplitude $\delta$}
\begin{figure}[tb]
\centerline{\epsfig{file=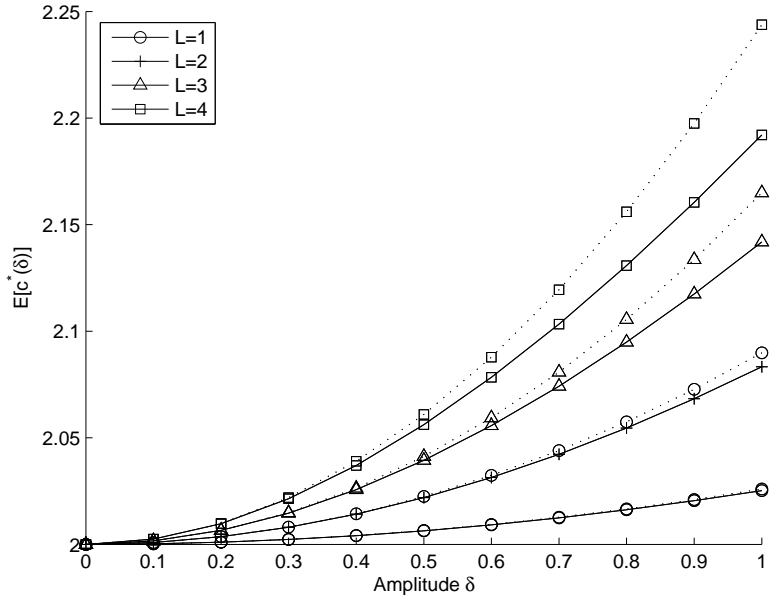,width=330pt}}
\caption{Average enhancement of minimal speed in small amplitude shears.
The solid curves are the numerically computed values;
the dashed curves are given by formula (\ref{e36}).}
\label{fig2}
\end{figure}

\begin{figure}[tb]
\centerline{\epsfig{file=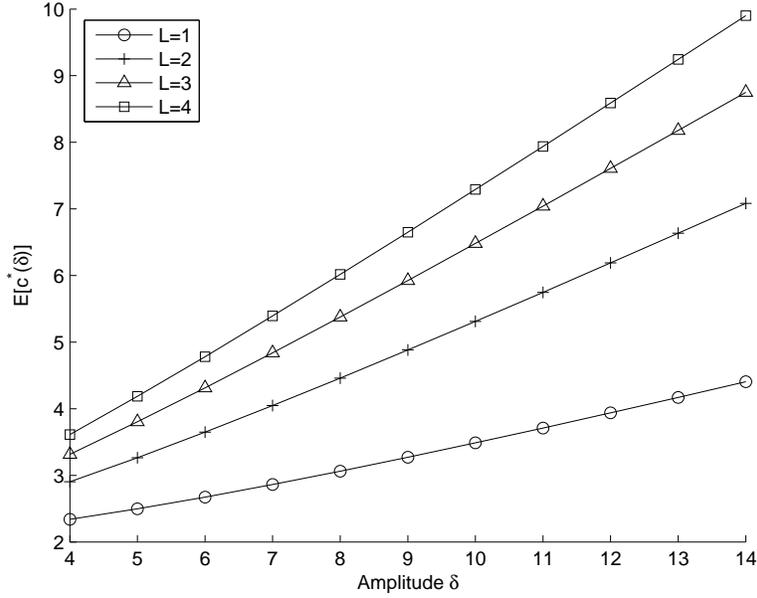,width=330pt}}
\caption{Average enhancement of minimal speed, large amplitude shears.}
\label{fig3}
\end{figure}

In Figure \ref{fig2} and Figure \ref{fig3}, we show the results of a simulation using $N=100,000$ realizations of a shear in small
and large root mean square (rms) amplitudes, respectively.
As shown later in Figure \ref{fig13}, we find that the choice of $N=100,000$ samples
was more than enough to obtain good convergence of the speed distribution functions.
The covariance function of the process is $E[b(y) b(s)] = 2e^{-4|y-s|}$. In each plot,
we show multiple curves, corresponding to various domain sizes. In Figure \ref{fig2},
corresponding to small $\delta$, the solid curves are the numerically computed values;
the dashed curves are given by formula (\ref{e36}).
We find that the enhancement of the minimal speed scales quadratically
for small amplitudes and linearly for large amplitudes. The computed exponents are shown in Table \ref{tab1}.

\begin{table}[htb]
\footnotesize
\caption{Computed scaling exponents $p$ for $E[c^*(\delta)] = c^*_0 + O(\delta^p)$.}\label{tab1}
\begin{center}
\begin{tabular}{| c | c | c  | c | c |}
\hline &  $L = 1.0 $ &  $L = 2.0 $ & $L = 3.0$ & $L = 4.0$  \\
\hline $\delta \ll 1$  &  $2.00 $ &  $1.98  $ & $1.96 $ & $1.93 $  \\
\hline $\delta \gg 1$  &  $1.09 $ &  $ 1.05 $ & $1.04 $ & $1.03 $  \\
\hline\end{tabular}
\end{center}
\end{table}

\begin{figure}[tb]
\centerline{\epsfig{file=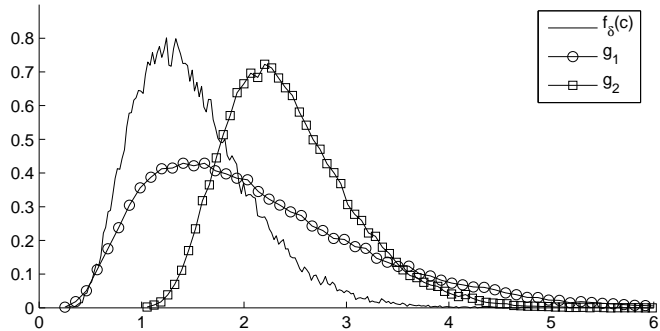,width=300pt}}
\caption{Distributions of speed enhancement, lower and upper bounds at $\delta = 50$.}
\label{fig15}
\end{figure}

As an application of the numerical results for large amplitudes, we compare the distribution of $c^*(\delta)$ with the distributions of the random variables $g_1(\omega) = 2 \kappa f'(0) + \delta \norm{b}_\infty$ and $g_2 (\omega) = 2 \sqrt{\kappa + \frac{\delta^2}{\kappa} \norm{\nabla \chi}_\infty}$, where
$\kappa $ is the diffusion constant (equal to 1 in equation (\ref{e1})).
Theorem 2 of \cite{HZ1} shows the upper bound:
\be
c^*(\delta,\omega) \leq \min ( g_1(\omega) , g_2(\omega))
\ee
and $g_1 < g_2$ provided that $\sqrt{\frac{\kappa}{f'(0)}}$ is sufficiently small, depending on the realization. In Figure \ref{fig15} we compare the distributions for $c^*(50)$, $g_1$, and $g_2$ for $\kappa = 0.01$. The asymmetry
is seen in all three curves.

%In equation (\ref{e37}), we see that for fixed $a>0$,
%\be
%\lim_{L \to \infty} E[\avg{\abs{\chi_x}^2}]  = +\infty.
%\ee
%So, by (\ref{e36}), we should expect that the enhancemement of the minimal speed increases linearly with the cylinder width at small amplitudes.

\subsubsection{Dependence on Covariance}
Next, we consider the effect of the covariance on the enhancement of the minimal speed. The covariance $E[b(y)b(s)]$ is a function of $\abs{t}=\abs{s-y}$, so we will write $V(t) = E[b(y)b(s)]$. By choosing $r = \sqrt{2} \alpha^{3/4}$,
we constructed O-U processes with covariances given by
\be
V(t) = \sqrt{\alpha} e^{-\alpha \, |t|} \quad.
\ee
By this choice of $r$, the $L^2$ norm of $V(t)$ remains constant as $\alpha$ changes, so that the total energy in the power spectrum of
the signal remains constant. Since $\frac{r^2}{2a} = \sqrt{a}$,
 we see from equation (\ref{e36}) that for fixed $L$,
\be
\lim_{\alpha \to +\infty} E[\avg{\abs{\chi_x}^2}] =
\lim_{\alpha \to 0^+} E[\avg{\abs{\chi_x}^2}] = 0
\ee
and that $E[\avg{\abs{\chi_x}^2}]$ achieves a maximum for some finite value of $\alpha \in (0,\infty)$.
This suggests that there is some optimal $\alpha$, depending on the domain size $L$, such that
the enhancement of $E[c^*(\delta)]$ is maximized.

Fixing the grid spacing $dx=0.002$, we computed the expected value $E[c^*(\delta)]$ for a
range of $\alpha$ and for $L= 1.0, 2.0, 3.0, 4.0$. Note that for each $\alpha$, we must
choose the initial points $b_0$ to have variance $E[b_0^2]= \sqrt{\alpha}$ so that the
process remains stationary for each $\alpha$.
Varying the covariance does not effect the order of the scaling in $\delta$. That is,
in each case the enhancement scales like $O(\delta^2)$ for small $\delta$ and $O(\delta)$ for large $\delta$,
as in the preceding simulation.

\begin{figure}[tb]
\centerline{\epsfig{file=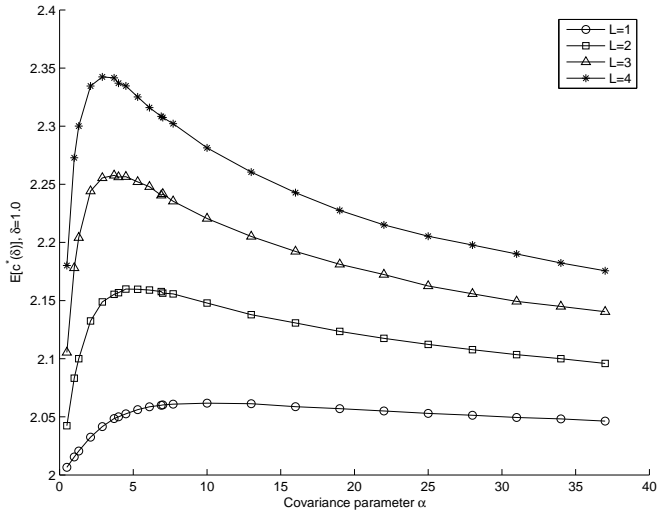,width=250pt}}
\caption{Effect of covariance on minimal speed enhancement at $\delta = 1.0$.}
\label{fig6a}
\end{figure}

\begin{figure}[tb]
\centerline{\epsfig{file=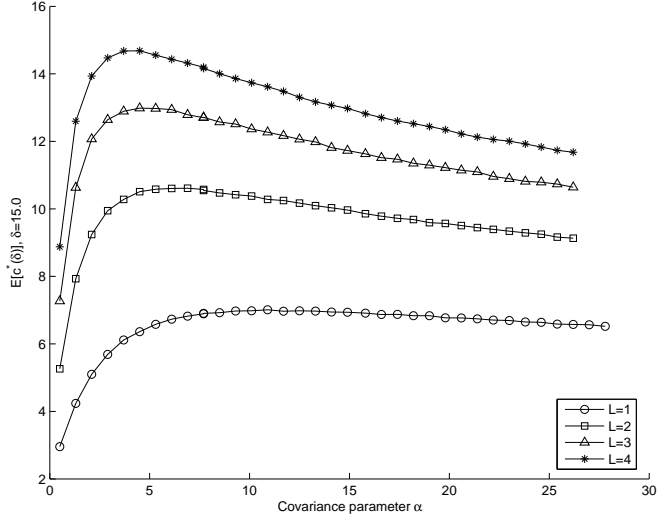,width=250pt}}
\caption{Effect of covariance on minimal speed enhancement at $\delta=15.0$.}
\label{fig8}
\end{figure}

Figure \ref{fig6a} shows the enhancement $E[c^*(\delta)]$ for a fixed $\delta=1.0$ and a range of $\alpha$.
Figure \ref{fig8} shows the results of the same computation for $\delta=15.0$,
corresponding to the large $\delta$ regime. In this case, formula  (\ref{e36})
is no longer valid. Nevertheless, we see the same effect as in the small
amplitude regime: the existence of an optimal covariance parameter $\alpha$.

This effect can be interpreted in terms of $V(t)$ and its Fourier transform or power spectrum:
\be
\hat V(w) = \sqrt{\frac{2}{\pi \alpha}} \left( 1 + \left(\frac{w}{\alpha}\right)^2 \right)^{-1}.
\ee
As $\alpha \to 0$, $\hat V(w)$ concentrates at the origin, and so the energy of the shear process
is concentrated more in the large scale spatial modes. The domain $\Omega$,
to which the process is restricted, is bounded, and variations over a length scale that
is much greater than the diameter of $\Omega$ have little effect on the average enhancement of the front.
As a result, $E[c^*]$ decreases as $\alpha \to 0^{+}$.  In the other limit
$\alpha \to \infty$, $\hat V$ spreads out so that the energy
over any finite band of frequencies goes to zero, causing $E[c^*]$ to decrease as well.
Note that $V \to 0$ in $L^1$ as $\alpha \to \infty$, so even though $\hat V$ spreads out more
uniformly as $\alpha \to 0$, the family of processes does not converge to white noise,
whose covariance function equals the Dirac delta function.

\subsubsection{Speed Distribution}

\begin{figure}[tb]
\centerline{\epsfig{file=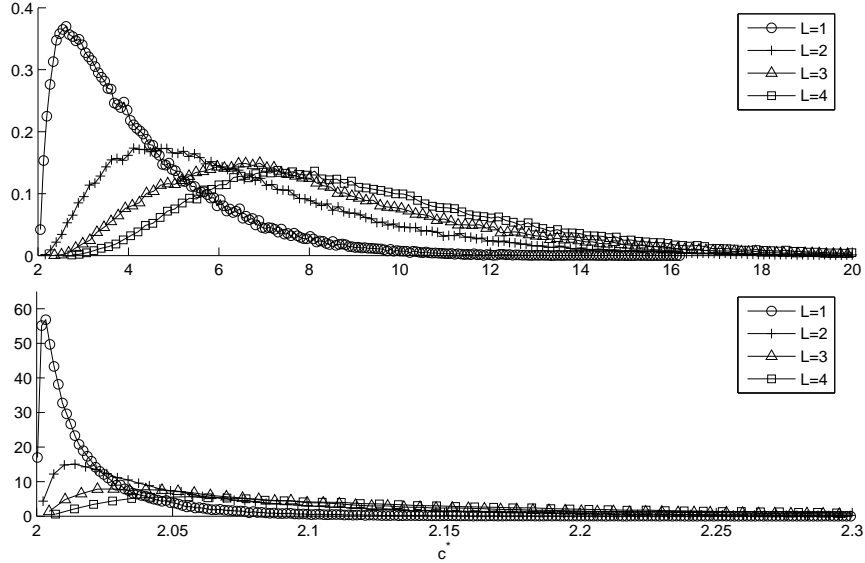,width=330pt}}
\caption{Computed probability distribution functions (pdf's) of enhancement $M(\delta)$, $\delta=1.0$ (bottom), $\delta=14.0$ (top).}
\label{fig4}
\end{figure}

For a fixed $\delta = 1$ and $\delta = 14$ (corresponding to small and large amplitudes), we computed the distributions of the numerically computed values $M(\delta)$. The distributions are shown in Figure \ref{fig4}. To compute these distributions, we partition the range of values into $Q$ disjoint intervals: $\left\{[x_j, x_{j+1})\right\}_{j=1}^Q$. Then, we let
\be
pdf(x) = \frac{1}{N} \sum_{i=1}^N \frac{\chi_j(M_i(\delta))}{ (x_{j+1} - x_j)} \quad \text{if} \;\; x \in [x_j,x_{j+1})
\ee
where $\chi_j(x)$ is the characteristic function of the interval
$[x_j,x_{j+1})$. The distributions in Figure \ref{fig4} were computed with $N=100,000$ samples, and $Q=300$.

\begin{figure}[tb]
\centerline{\epsfig{file=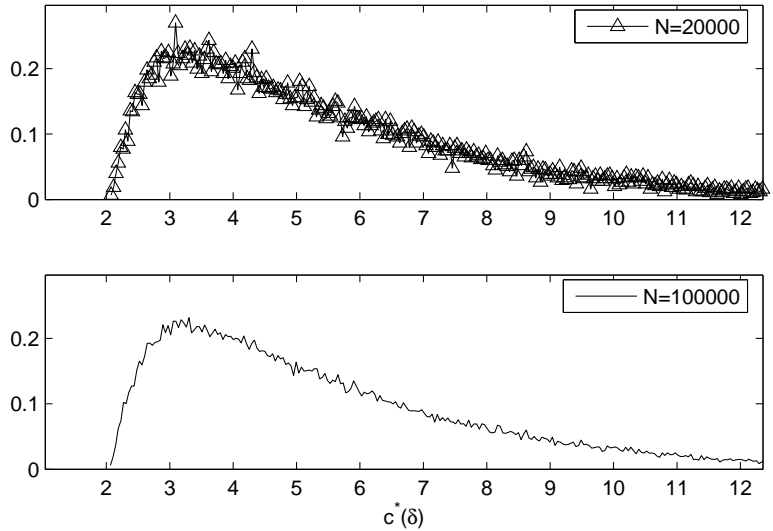,width=350pt}}
\caption{Convergence of speed enhancement distribution at $\delta = 14.0$. }
\label{fig13}
\end{figure}

The values $M(\delta)$ are the enhancement of the minimal speeds due to the variation of the shear, after the effect of the mean field has been subtracted off.  Since a mean zero shear always enhances the minimal speeds, we should expect $M(\delta) > 0$ for all $\delta$, for all realizations. As already noted, Figure \ref{fig13} shows good convergence of the distributions when using $N=100,000$ samples.

\subsubsection{Comparison with Direct Simulation}
When comparing the results of the previous sections with results from direct simulation, we find that computing the minimal speeds using the
variational formula offers significant advantages. As in \cite{NX1}, we
first computed the quantities $E[c^*(\delta)]$ via direct simulation of the original
equation (\ref{e1}) on a truncated domain using an explicit second-order upwind finite difference scheme. We chose the diffusion constant in the direct simulations to be $\kappa = 0.025$ and the grid spacing to be $dx = 0.05$, $dt = 0.004$. For each realization of the shear, we evolved the solution
of (\ref{e1}) for a sufficiently long time until the front moves at a more or less constant speed (see \cite{NX1} for details of the method). In this way, we approximated the
minimal speed for each realization; then, we repeated the process for a large ensemble of shears to compute the expectations $E[c^*(\delta)]$.

\begin{figure}[tb]
\centerline{\epsfig{file=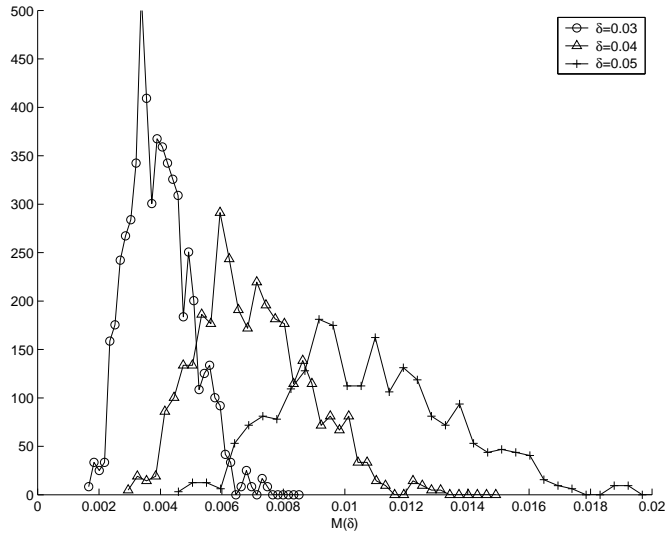,width=250pt}}
\caption{Distribution of speed enhancement via direct simulation, KPP nonlinearity, $\delta = 0.5$.}
\label{fig9}
\end{figure}

\begin{figure}[tb]
\centerline{\epsfig{file=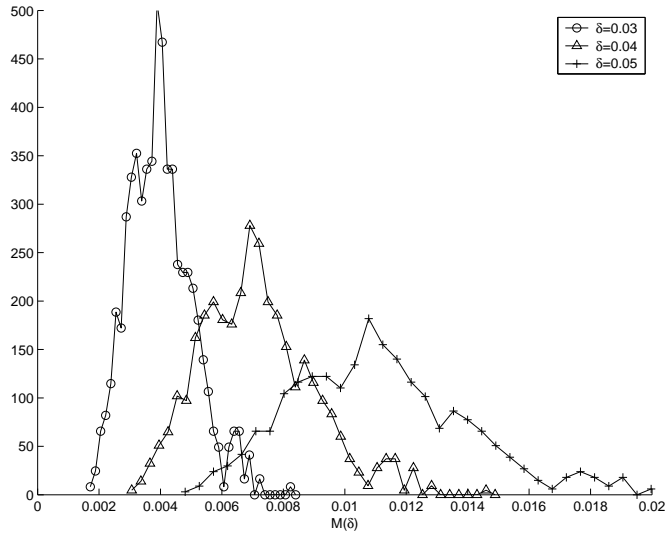,width=250pt}}
\caption{Distribution of speed enhancement via direct simulation, combustion nonlinearity, $\delta = 0.5$.}
\label{fig10}
\end{figure}

\begin{figure}[tb]
\centerline{\epsfig{file=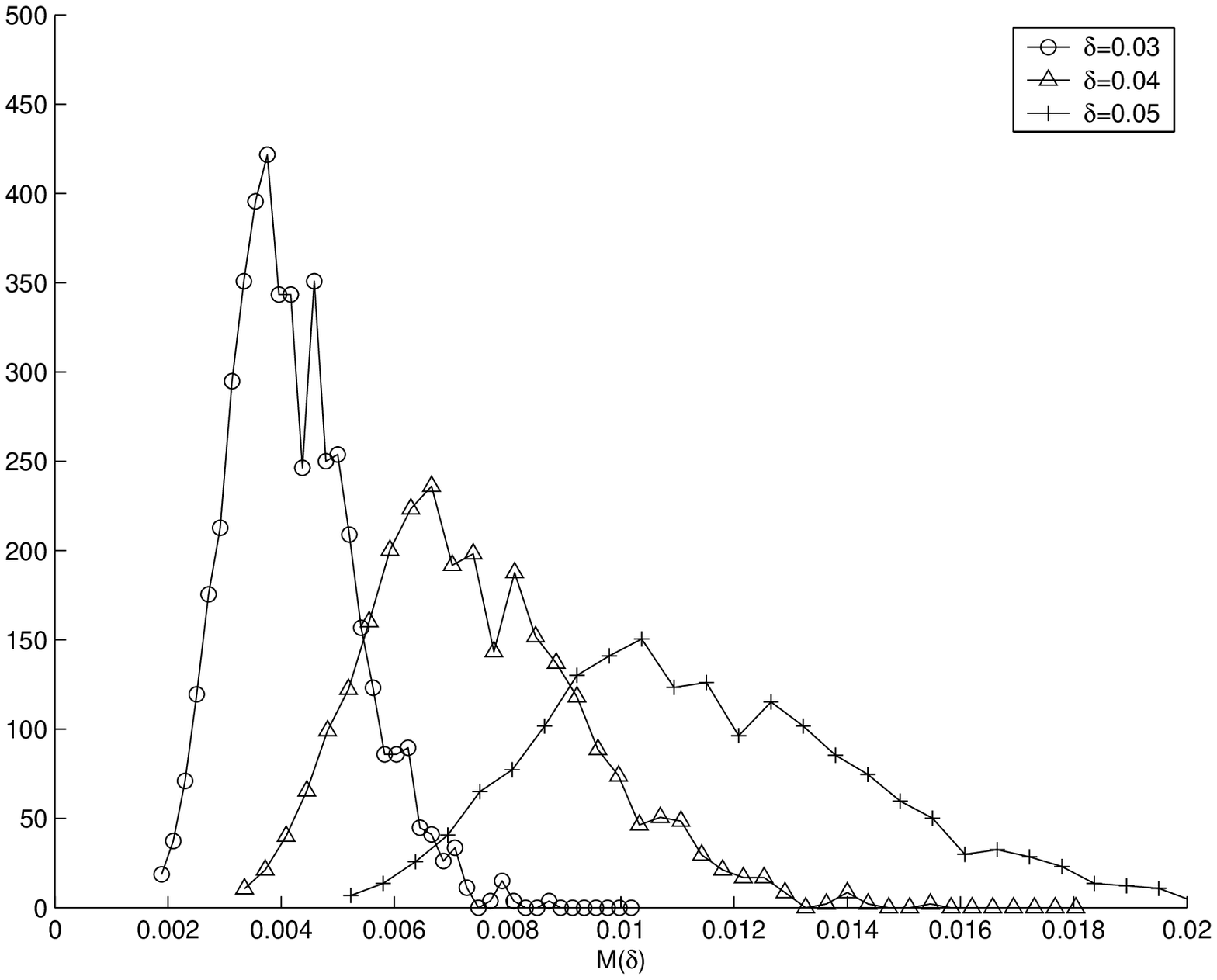,width=250pt}}
\caption{Distribution of speed enhancement via direct simulation, bistable nonlinearity, $\delta = 0.5$.}
\label{fig11}
\end{figure}

We considered three nonolinearities in our direct simulations:
\begin{itemize}
\item KPP nonlinearity: $f(u) = u(1-u)$.
\item Combustion nonlinearity: $f(u) \equiv 0$ for $u \in [0,\theta]$, and $f(u) > 0$ for $u \in (\theta, 1]$, for some $\theta \in (0,1)$.  Also, $f'(1) < 0$.
\item Bistable nonlinearity: $f(u) = u(1-u)(u-\mu)$ for some $\mu \in (0,1/2)$.
\end{itemize}
Although there is no known variational formula as simple as formula (\ref{e8}) for the bistable and combustion nonlinearities, the expansion (\ref{e16}) holds for
each of the nonlinearites (see also Theorem 4.2 of \cite{HPS}).  Therefore, we should expect that for small $\delta$,
the computed values
\be
\frac {M(\delta,\omega)}{c^*_0} \approx \frac{\delta^2}{2 \abs{\Omega}} \int_\Omega \abs{\nabla \chi(y,\omega)}^2\;dy + O(\delta^3)
\ee
have approximately the same distribution, independent of the nonlinearity.
For each nonlinearity, we computed more than 700 realizations of the shear and
evolved the solution, as described in \cite{NX1}.  A
 relatively smaller number of samples
was due to the long time required to compute each
sample (not a problem with the variational formula in the KPP case).
After subtracting off the linear part of the enhancement (due to mean field), we computed the
enhancement exponent as described in Section 3.  The results of the direct simulation concur that, for small amplitudes, $E[c^*(\delta)]$ scales like $O(\delta^2)$.
Figures \ref{fig9}, \ref{fig10}, and \ref{fig11} show the distributions of the computed values $\frac{M(\delta)}{c^*_0}$ for the KPP, combustion, and bistable nonlinearities, respectively, for $L=1$.
We see that the distributions are roughly the same, independent of the nonlinearity, as should be expected.
While computationally expensive, the direct simulation method reveals the universal $O(\delta^2)$ scaling
of the front speeds for small $\delta$ (albeit on smaller ensembles than in the variational computations),
for each of the nonlinearities considered.

Because the shears are random and may greatly distort the wave front, one major
challenge in accurately approximating the (minimal) speeds is tracking
the widely varying front region over a long time.  In the region of the front, gradients are relatively large, so accurately tracking the front requires either a very fine
uniform grid spanning a large domain or some kind of adaptive-mesh scheme.  In either case, accurate direct simulation
is prohibitively expensive compared to the simple variational method.  In contrast to the direct simulation method, the variational formula allowed us to
compute a much larger number of samples in a fraction of the time and avoid
the inaccuracies resulting from truncation of the channel doman.

Finally, we note that if $\Omega$ is unbounded, then the quadratic
asymptotic behavior of $c^*(\delta)$ cannot hold in general. For example, in \cite{Xin2},
it was shown that if the channel $R \times [0,L]$ is replaced by $R^2$ then the
front velocities obtained through a comparable variational principle diverge
due to the almost sure growth of the running maximum of the process $b(y,\omega)$.
This effect can be clearly seen in our results, as $E[\avg{\abs{\nabla \chi}^2}]$
given by (\ref{e36}) diverges as $L \to +\infty$.
In our application of Theorem 2.1 to a stationary Gaussian process, the
boundedness of the domain $\Omega$ is crucial in order to achieve the
necessary bound on $\norm{b}_\infty$ (through the martingale moment inequality),
and in turn the bound on $c^*$.

%The hypotheses of Theorem 2.1 can
%be achieved for processes that are nowhere differentiable (fractal).
%we do not observe the anomalous speed enhancement described in \cite{MS2} for fronts
%propagating in $R^2$ with turbulent velocity fields.

\section{Conclusions}
Sufficient moment conditions are obtained to ensure the quadratic (linear) KPP average front speed enhancement
through small (large) rms random shear flows in channel domains. The conditions are realized by the Ornstein-Uhlenbeck
process for which an explicit average front speed formula is derived. The variational principle based
computation is carried out for the speed ensemble. Numerically computed speed enhancement is in
agreement with theoretical analysis, and provides data for studying speed distributions and dependence on
shear covariance. Comparison with direct simulations of random fronts in case of non-KPP nonlinearities
showed the same enhancement scaling laws. It is interesting to investigate front speeds through
time dependent random shear flows or non-shear random flows in channel domains for future studies.

\vspace{.2 in}

\section{Acknowledgements}
The work was partially supported by NSF grant ITR-0219004.
J. X. would like to acknowledge a fellowship from the John Simon Guggenheim Memorial Foundation, and a
Faculty Research Assignment Award at UT Austin.
J. N. is grateful for support through a VIGRE graduate fellowship at UT Austin.
Both authors wish to thank J. Wehr for helpful communications.

%\newpage
\vspace{.2 in}

\bibliographystyle{plain}

\end{document}